\documentclass[11pt]{article}
\usepackage{amsmath}
\usepackage{graphicx}
\usepackage{amsfonts}
\usepackage{amssymb}
\def\LARGE{\Large}
\def\Large{\large}

\begin{document}
\begin{center}
{\LARGE{\bf {\Large Allee effects in a Ricker-type predator-prey
system }}}
\end{center}

\begin{center}Yunshyong Chow$^1$ and Sophia R.-J.  Jang$^2$\end{center}

\noindent 1. Institute of Mathematics, Academia Sinica, Taipei
10617, Taiwan

\noindent 2. Department of Mathematics and Statistics, Texas Tech
University, Lubbock, TX 79409-1042 U.S.A.

\vspace{0.1in}

\noindent\textbf{Abstract.} We study a discrete host-parasitoid
system where the host population follows the classical Ricker
functional form and is also subject to Allee effects. We determine
basins of attraction of the local attractors of the single
population model when the host intrinsic growth rate is not large.
In this situation, existence and local stability of the interior
steady states for the host-parasitoid interaction are completely
analyzed. If the host's intrinsic growth rate is large, then the
interaction may support multiple interior steady states. Linear
stability of these steady states  is provided.

 \vspace{0.1in}


\vspace{0.1in}


\section{Introduction}
\setcounter{equation}{0}

The Allee effect, refers to the reduced fitness or the decline in
population growth at low population densities or sizes, was first
observed by Allee \cite{allee}.  It has a significant impact on
population survival when the population is at low level. There has
been a rebound of interest in Allee effects recently due to
fragmentation of habitats, invasions of exotic species, biological
control of pest, etc., all involved with small populations. The
treatise by Courchamp et al. \cite{courchamp} provides interesting
biological motivation and many mathematical models of Allee
effects.

The Ricker model \cite{ricker} is one of the classical equations
used to model fish populations and the function has been applied
extensively to model other populations as well. Its dynamics are
complicated. As the intrinsic growth rate of the population is
increased, the map undergoes a cascade of period-doubling
bifurcations and eventually becomes chaotic. Due to its
importance, the Ricker growth function has been incorporated with
Allee effects to study population dynamics by several researchers.
Li et al. \cite{li} propose a single population model with Allee
effects by assuming a Holling type II functional form for the
birth or growth of the population, termed the extended Ricker
model. The authors \cite{li}  compare the dynamics of the extended
Ricker model with the classical Ricker model and conclude that
Allee effects have the stabilizing effects on population dynamics.
In \cite{elaydi1}, Elaydi and Sacker study a modified Ricker map
that exhibits Allee effects. The modified Ricker map is then
generalized to include inter-specific interactions of two
populations  by Livadiotis and Elaydi \cite{liva}.  Jang and
Diamond \cite{jang} on the other hand study a predator-prey system
with Allee effects occurring in the prey that grows according to
the Ricker function. More recently, Kang proposes a two species
competition model with Allee effects occurring in both populations
\cite{kang}, where each single population in the absence of the
other population follows the growth of the Ricker model and is
also subject to Allee effects. There are other mathematical models
of Allee effects, including \cite{cushing, dennis, jang1,
schreiber, thieme, zhou, zhou1} and references cited therein.

In each of these Ricker-type models \cite{elaydi1, jang,  kang,
li, liva}, however, the mechanism of Allee effects is modeled
differently from the present single population  model studied in
this manuscript. Specifically, the way that Allee effects are
incorporated into the Ricker equation is the same as that of the
deterministic single population model proposed by Allen et al.
\cite{allen1} and by Liebhold and Bascompte \cite{liebhold}. To
the best of our knowledge, the resulting scalar equation with
Allee effects has not been analyzed previously. We derive a
threshold $r_0$ in terms of the Allee threshold such that the
single population possesses only equilibrium dynamics if its
intrinsic growth rate $r$ is smaller than the threshold $r_0$.
Basins of attraction of the local attractors are then explicitly
determined. The population exhibits a period-doubling bifurcation
at the critical threshold $r=r_0$. Based on this single population
equation, we build a host-parasitoid system. We study existence
and local stability of the interior steady states. Under certain
parameter regime, global dynamics are provided for this
two-dimensional system. Numerical simulations are also performed
using Matlab to study the model.

The host-parasitoid model is presented in Section 2 and the single
host population is analyzed in Section 2.1. Section 3 provides
analysis of the two-dimensional host-parasiotid system. In
particular, preliminary results are given in Section 3.1 and
Sections 3.2 and 3.3 deal with the cases of $r<r_0$ and $r>r_0$
respectively. The final section provides a brief summary.

\section{The model and the host equation}
\setcounter{equation}{0} Let $x(t)$ and $y(t)$ be the host and
parasitoid populations at time $t=0, 1, \cdots$, respectively.
Parameters  $r$ and $K$ are the intrinsic growth rate and the
carrying capacity of the host population, respectively. Inspired
by the  studies in \cite{amar, keitt}, where the continuous-time
logistic model $N^\prime(t)=rN(1-N/K)$ is modified to yield
$N^\prime(t)=rN(1-N/K)(N-a)$ for the mechanism of Allee effects,
Liebhold and Bascompte \cite{liebhold} derive the discrete-time
version of this type of model for organisms with non-overlapping
generations via
$$\ln(\tilde\gamma(t))=r(1-x(t)/K)(x(t)-a),$$ where $\tilde\gamma(t)=x(t+1)/x(t)$
is the change in population density. As a result, one obtains
$$x(t+1)=x(t)e^{r(1-x(t)/K)(x(t)-a)}$$ for the single host
population that is subject to Allee effect with the Allee
threshold given by $a$, $0<a<K$.
 The per capita growth rate of the host population is less than
$1$ if its population level is below $a$ and this growth rate
exceeds $1$ if the population level is between $a$ and the
carrying capacity $K$. The per capita growth rate is smaller than
$1$ again if the population level is larger than $K$ due to
intra-specific competition.

The parasitoid population feeds on the host population
exclusively. Parameter $\beta>0$ is the conversion of prey to
predator and $b>0$ represents parasitoid's searching and handling
efficiency on the host. The host-parasitoid interaction is given
by the following system:
\begin{equation} \left \{
\begin{array}{l}
x(t+1)=x(t)e^{r(1-x(t)/K)(x(t)-a)}e^{-by(t)}\\[1ex]
y(t+1)=\beta x(t)(1-e^{-by(t)}),\\
 \end{array}\right.
\end{equation}
where $0<a<K$ and $r, b, \beta>0$. There are $5$ parameters in
model (2.1). We  can rescale the state variables and parameters to
eliminate some parameters. Let
$$\hat x(t)=x(t)/K, \ \hat y(t)=by(t), \ \hat r=rK, \
\hat a=a/K, \ \hat \beta=\beta bK.$$ Then system (2.1) is
converted to
\begin{equation*}
\left \{
\begin{array}{l}
\hat x(t+1)=\hat x(t)e^{\hat r(1-
\hat x(t))(\hat x(t)-\hat a)}e^{-\hat y(t)}\\[1ex]
\hat y(t+1)=\hat\beta \hat x(t)(1-e^{-\hat y(t)}),
 \end{array}\right.
\end{equation*} where $0<\hat a<1$. By ignoring the hats,
we obtain the following system with only $3$ parameters
\begin{equation}
\left \{
\begin{array}{l}
x(t+1)=x(t)e^{r(1-x(t))(x(t)-a)}e^{-y(t)}\\[1ex]
y(t+1)=\beta x(t)(1-e^{-y(t)}),
 \end{array}\right.
\end{equation} where $r, \beta>0$ and $0<a<1$, and  initial conditions are nonnegative. The parameters have different biological meanings than those in model (2.1).
We shall study system (2.2) in section 3 and we will focus on the
analysis of single host population model in section 2.1.

\subsection{The host population} If $y(0)=0$  in (2.2), we have $y(t)=0$ for
$t\geq 1$  and (2.2) reduces to a one-dimensional equation
\begin{equation}
x(t+1)=x(t)e^{r(1-x(t))(x(t)-a)}.
\end{equation}
Model (2.3) has three steady  states: $0, \ a$ and $1$. Let $f(x)$
denote the map of (2.3),
\begin{equation} f(x)=xe^{r(1-x)(x-a)}.\end{equation} Then
\begin{equation}f^\prime(x)=e^{r(1-x)(x-a)}\left(1+xr(1-2x+a)\right),\end{equation}
and
\begin{equation}
f^{\prime\prime}(x)=r\left(rx(1+a-2x)^2+3(1+a-2x)-(1+a)\right)e^{r(1-x)(x-a)}.
\end{equation}
Since $0<f^\prime(0)=e^{-ra}<1$, steady state $0$ is always
locally asymptotically stable. On the other hand,
$f^\prime(a)=1+ar(1-a)>1$, and hence steady state $a$ is a
unstable. Since $f^\prime(1)=1+r(a-1)<1$,  $1$ is locally
asymptotically stable if $r<\displaystyle\frac{2}{1-a}$ and
unstable if $r>\displaystyle\frac{2}{1-a}$. We use $r$ as the
bifurcation parameter to study (2.3).

Denote $r_0$ in terms of the Allee threshold $a$:
\begin{equation} r_0=\displaystyle\frac{2}{1-a}.\end{equation}
Then $f^\prime(1)|_{r=r_0}=-1$ and it is expected that a
period-doubling bifurcation occurs at $r=r_0$. In the following,
we first verify that steady state $1$ is locally asymptotically
stable when $r=r_0$. Consider the Schwarzian derivative of $f$,
$$Sf(x)=\displaystyle\frac{f^{\prime\prime\prime}(x)}{f^\prime(x)}-\displaystyle\frac{3}{2}
\left(\displaystyle\frac{f^{\prime\prime}(x)}{f^\prime(x)}\right)^2.$$
By \cite[Theorem 1.16]{elaydi}, steady state $1$ is locally
asymptotically stable if $Sf(1)<0$ and  unstable if $Sf(1)>0$. A
direct computation yields

\begin{equation}f^{\prime\prime}(1)|_{r=r_0}=-2r_0 \mbox{ and }
f^{\prime\prime\prime}(1)|_{r=r_0}=8r_0-2ar_0.\end{equation}
Therefore $Sf(1)|_{r=r_0}=2r_0(a-4-3r_0)<0 \mbox{ for } 0<a<1$ and
steady state $1$ is locally asymptotically stable when $r=r_0$. To
prove period-doubling bifurcation, we
 apply \cite{elaydi} to verify that
\begin{equation}\displaystyle\frac{\partial^2f}{\partial x\partial
r}(1,r_0)\neq 0 \mbox{ and }
\displaystyle\frac{1}{2}\left(\displaystyle\frac{\partial^2
f}{\partial
x^2}(1,r_0)\right)^2+\displaystyle\frac{1}{3}\displaystyle\frac{\partial^3f}{\partial
x^3}(1,r_0)\neq 0.\end{equation} Notice
$\displaystyle\frac{\partial^2f}{\partial x\partial
r}(1,r_0)=a-1<0$ and  the first expression in (2.9) is satisfied.
Substituting (2.8) into the second expression of (2.9) arrives at
$\displaystyle\frac{2r_0}{3}(3r_0+4-a)>0$ and therefore (2.9)
holds. As a result, a period-doubling bifurcation occurs at
$r=r_0$. We now summarize.

\vspace{0.1in}

\noindent{\bf Proposition 2.1} {\em Model (2.3) has three steady
states $0, \ a$ and $1$, where $0$ is locally asymptotically
stable and $a$ is unstable. Steady state $1$ is locally
asymptotically stable if $r\leq r_0$ and unstable if $r>r_0$.
Moreover, equation (2.3) undergoes a period-doubling bifurcation
at $r=r_0$.}

\vspace{0.1in}

It follows from Proposition 2.1 that model (2.3) has $2$-cycle
solutions when $r>r_0$ and is close to $r_0$. To further study the
existence of $2$-cycles, we consider the second composition of
$f$, $f^2(x):=f(f(x))$, where \begin{equation} f^2(x)=x
e^{\displaystyle r(1-x)(x-a)}e^{\displaystyle
r(1-xe^{r(1-x)(x-a)})(xe^{r(1-x)(x-a)}-a)}.\end{equation} Two-
cycle solutions are fixed points of the map $f^2$ that are not
fixed points of $f$. Consequently, components $x$ of any
$2$-cycles satisfy
\begin{equation} (1-x)(x-a)+
\left(1-xe^{r(1-x)(x-a)}\right)\left(xe^{r(1-x)(x-a)}-a\right)=0.\end{equation}
Clearly $a$ and $1$ are solutions of (2.11) and we look for other
positive solutions. Let $\bar g(x)$ denote the left hand side of
(2.11). Then $\bar g(0)=-2a<0$, $\bar g(a)=\bar g(1)=0$,
$\displaystyle\lim_{x\rightarrow\infty}\bar g(x)=-\infty$ and
\begin{equation}
\bar
g^\prime(x)=a-2x+1-e^{r(1-x)(x-a)}[1+rx(1-2x+a)][2xe^{r(1-x)(x-a)}-a-1].\end{equation}
Notice $\bar g^\prime(a)=(1-a)(2+ra(1-a))>0$, and $\bar
g^\prime(1)=(a-1)(2+r(a-1))>0$ if and only if $r>r_0$. We conclude
that (2.3) has a positive $2$-cycle ${\cal S}=\{\bar x_1, \bar
x_2\}$ if $r>r_0$, where $a<\bar x_1<1<\bar x_2$.

Notice that $f$ attains its maximum at
\begin{equation}x_m=\displaystyle\frac{1+a}{4}+\displaystyle\frac{\sqrt{(a+1)^2r^2+8r}}{4r},
\end{equation}
where $x_m>\displaystyle\frac{1+a}{2}>a$.  Moreover, $x_m>1$ if
and only if $r<r_0/2$, and $x_m$ is a decreasing function of $r$.
See Figure 1(a)--(b) for the illustration.  Notice that there
exists a unique $x_a>\mbox{max}\{1, \ x_m\}$ such that $f(x_a)=a$.
 Let ${\cal B}(x)$ denote the basin of attraction
of the steady state $x$. If $x(0)<a$, then $x(1)=f(x(0))<f(a)=a$
and $x(1)=f(x(0))<x(0)$. Inductively, $\{x(t)\}_{t=0}^{\infty}$ is
a decreasing sequence and bounded below by $0$. Thus,
$\displaystyle\lim_{t\rightarrow\infty}x(t)=0$ by the continuity
of $f$, i.e., $[0, a)\subset{\cal B}(0)$. Furthermore, since $f$
is decreasing on $(x_m, \infty)$, we have $f(x)<a$ if $x>x_a$.
Therefore, $\displaystyle\lim_{t\rightarrow\infty}x(t)=0$ if
$x(0)>x_a$ and $[0, a)\cup(x_a, \infty)\subset {\cal B}(0)$ for
any $r>0$.

In the following, we   prove that solutions of (2.3) either
converge to $0, a $ or $1$ when $r<r_0$. We separate the analysis
into two cases: $r\leq r_0/2$ and $r_0/2<r<r_0$. Notice that
$1<x_m$ and $1>x_m$ when $r<r_0/2$ and $r>r_0/2$, respectively.

\vspace{0.1in}

\noindent{\bf Theorem 2.2} {\em If $r\leq r_0/2$, then  ${\cal
B}(0)= [0, a)\cup(x_a, \infty)$ and ${\cal B}(1)= (a, x_a)$.}

\vspace{0.1in}

\noindent{\bf Proof.} First assume   $r<r_0/2$. Then $x_m>1$. If
$x\in(a, 1)$, then $f(x)>x>a$ and $f(x)<1$ and hence
$\displaystyle\lim_{t\rightarrow\infty}f^t(x)=1$.  If $x\in (1,
x_m]$, then $f(x)<x\leq x_m$ and $f(x)>1$ hold and thus
$\displaystyle\lim_{t\rightarrow\infty}f^t(x)=1$. If $x\in (x_m,
x_a)$, then $x_m>f(x_m)>f(x)>f(x_a)=a$, i.e., $f(x)\in (a, x_m)$,
and hence $\displaystyle\lim_{t\rightarrow\infty}f^t(x)=1$.
Similar analysis can be applied to the case when $r=r_0/2$ and the
proof is complete. \rule{2mm}{2mm}

\vspace{0.1in}

Suppose now $r_0/2<r<r_0$. Then $x_m<1$ and $f^\prime(1)<0$.
Recall that $\displaystyle\lim_{t\rightarrow\infty}x(t)=0$ if
$x(0)\in[0, a)\cup (x_a,\infty)$. We prove that $f(x_m)\geq x_a$
cannot occur and the basin of attraction of $1$ is $(a, x_a)$.

\vspace{0.1in}

 \noindent{\bf Theorem 2.3} {\em If $r_0/2<r<r_0$, then
 $f(x_m)<x_a$.  Moreover, ${\cal
B}(0)= [0, a)\cup(x_a, \infty)$ and ${\cal B}(1)= (a, x_a)$. }

\vspace{0.1in}

\noindent{\bf Proof.} Assume first that $f(x_m)<x_a$ holds. Then
$(a, x_a)$ is positively invariant for $f$, $f$ has a maximum at
$x_m\in (a,1)$ and $f$ is decreasing on $(x_m,\infty)$. We apply
Theorem 2.9(b) of \cite[page 50]{allen}. Specifically, since $f$
satisfies (i)-(vi) given in \cite[Theorem 2.9(b)]{allen}, $1$ is
globally asymptotically stable on $(a, x_a)$ if
\begin{equation} f(f(x))>x \mbox{ for } x\in[x_m, 1).
\end{equation} Using (2.10) and (2.11), (2.14) is equivalent to
\begin{equation} (1-x)(x-a)+(1-f(x))(f(x)-a)>0 \mbox{ for } x\in[x_m, 1).
\end{equation}
For $x\in(a, 1)$, let \begin{equation}
y=\displaystyle\frac{1+a+\sqrt{(1-a)^2+4(1-x)(x-a)}}{2}\end{equation}
be the solution of
\begin{equation} (1-x)(x-a)+(1-y)(y-a)=0
\end{equation} that is greater than $1$. Then (2.15) is equivalent to
\begin{equation} G(x):=y-f(x)>0 \mbox{ on } [x_m, 1),
\end{equation}
where $x_m$ depends on $r$ and $a$. As a result, we will need to
show that $G(x)>0 $ for $x\in [x_m, 1)$, $0<a<1$ and
$r_0/2<r<r_0$. Since $f(x)$ is strictly increasing in $r$ for
$x\in(a, 1)$, it is enough to show that (2.18) holds at $r=r_0$
and $0<a<1$. On the other hand, $G(x_m)>0$ is equivalent
$f(f(x_m))>x_m$  and $f(f(x_m))>x_m$ implies $f(f(x_m))>a$. This
later inequality is also equivalent to $f(x_m)<x_a$. Therefore, if
(2.18) is verified at $r=r_0$ and for all $a\in(0, 1)$, then
$f(x_m)<x_a$ holds and the earlier assumption can be dropped.
Consequently $1$ is globally asymptotically stable on $(a, x_a)$
is also proved. Therefore, we are  in a position to verify (2.18)
for $r=r_0$ and $0<a<1$. If we can show that $G(x)$ is strictly
decreasing on $[x_m, 1)$, then since $G(1)=0$ we have $G(x)>0$ on
$[x_m,1)$ and (2.18) is shown. Therefore, in the following we will
prove that $G(x)$ is strictly decreasing on $[x_m, 1)$ when
$r=r_0$ and $0<a<1$.

To this end, differentiating $y$ with respect to $x$, yields
\begin{equation}
y^\prime(x)=\displaystyle\frac{1+a-2x}{\sqrt{(1-a)^2+4(1-x)(x-a)}}<0
\mbox{ on } [x_m, 1),
\end{equation}
and
\begin{equation}
y^{\prime\prime}(x)=\displaystyle\frac{-4(1-a)^2}{\left((1-a)^2+4(1-x)(x-a)\right)^{3/2}}<0
\mbox{ on } [x_m, 1).
\end{equation}
Moreover, it is easy to see from (2.20) that $y^{\prime\prime}(x)$
is strictly decreasing on $[x_m, 1)$. Let \begin{equation}
P_3(x)=rx(1+a-2x)^2+3(1+a-2x)-(1+a).\end{equation} Then
$f^{\prime\prime}(x)=rP_3(x)e^{r(1-x)(x-a)}$ by (2.6) and
$$P_3(0)>0, \ P_3(\displaystyle\frac{1+a}{2})=-(1+a)<0,  \mbox{
and } P_3(1)=-2<0.$$ Since $$\displaystyle\lim_{x\rightarrow
-\infty}P_3(x)=-\infty \mbox{ and }
\displaystyle\lim_{x\rightarrow \infty}P_3(x)=\infty,$$ it follows
that $P_3(x)<0$ on $[x_m,1)$ and hence
\begin{equation}f^{\prime\prime}(x)<0 \mbox{ on }[x_m, 1).\end{equation}
A straightforward computation shows that
\begin{equation}
f^{\prime\prime\prime}(x)=rP_4(x)e^{r(1-x)(x-a)},
\end{equation}
where
\begin{equation}
P_4(x)=-\displaystyle\frac{r^2}{2}z^4-\displaystyle\frac{r^2(1+a)}{2}z^3+6rz^2+3r(1+a)z-6
\end{equation} and $z=2x-1-a$.
If we can show that $P_4(x)>0$ on $[x_m,1)$, then
$f^{\prime\prime\prime}(x)>0$ on $[x_m,1)$, which implies
$f^{\prime\prime}(x)$ is strictly increasing on $[x_m, 1)$. Since
$y^{\prime\prime}(x)$ is strictly decreasing on $[x_m ,1)$ by
(2.20), we see that
$G^{\prime\prime}(x)=y^{\prime\prime}(x)-f^{\prime\prime}(x)$ is
strictly decreasing on $[x_m, 1)$. Now,
$y^{\prime\prime}(1)=-\displaystyle\frac{4}{1-a}$ and
$f^{\prime\prime}(1)=-\displaystyle\frac{4}{1-a}$ at $r=r_0$ imply
 $G^{\prime\prime}(1)=0$. Hence $G^{\prime\prime}(x)>0$ on $[x_m,
1)$ and $G^\prime(x)$ is strictly increasing on $[x_m, 1)$. Since
$G^{\prime}(1)=y^\prime(1)-f^\prime(1)=0$, we must have
$G^\prime(x)<0$ on $[x_m, 1)$. As $G(1)=y(1)-f(1)=0$, we conclude
that $G(x)>0$ on $[x_m, 1)$.

It remain to prove that $P_4(x)>0$ on $[x_m, 1)$. Differentiating
$P_4(x)$ yields
$$
P_4^\prime(x)=-4r^2(2x-1-a)^2-3r^2(1+a)(2x-1-a)^2+24
r(2x-1-a)+6r(1+a)$$ and
$$P_4^{\prime\prime}(x)=-12r\left(2r(2x-1-a)^2+r(1+a)(2x-1-a)-4\right).$$
Solving $P_4^{\prime\prime}(x)=0$, we obtain
$$x_{\pm}=\displaystyle\frac{1+a}{2}-\displaystyle\frac{1+a}{8}\pm\displaystyle\frac{\sqrt{r^2(1+a)^2+32r}}{8r}.$$
Notice $$x_-<\displaystyle\frac{1+a}{2}<x_+, \mbox{ and } x_+<1
\mbox{ if and only if } 8(a-1)^2>0.$$ This later inequality holds
trivially since $a<1$ and hence $x_+<1$ is valid. Furthermore, the
graph of $y=P_4^{\prime\prime}(x)$ is concave down, and
$P_4^{\prime\prime}(x)$ is positive on $(x_-, x_+)$ and negative
on $(x_+, \infty)\cup(-\infty, \ x_-)$. As a result, the graph of
$y=P_4^\prime(x)$ has  critical points at $x=x_\pm$ and
$P_4^\prime$ is increasing on $(x_-, x_+)$ and decreasing on
$(x_+, \infty)\cup(-\infty, \ x_-)$. Since
$$P_4^\prime(\displaystyle\frac{1+a}{2})=6r(1+a)>0$$  and
$$P_4^\prime(1)=-16(1-a)-12(1+a)+48+\displaystyle\frac{12}{1-a}(1+a)>0,$$
we conclude that $P_4^\prime(x)>0$ on
$[\displaystyle\frac{1+a}{2}, 1]$. We need to determine the sign
of  $P_4(x_m)$. Since $x_m$ satisfies $1+rx(1+a-2x)=0$, which
implies $2xr(2x-1-a)-2=0$ and hence $x_m$ solves
\begin{equation} r(2x-1-a)^2+r(1+a)(2x-1-a)-2=0.\end{equation}
Using (2.25) and  long division, we obtain
\begin{equation} P_4(x_m)=-2r(1+a)(2x_m-1-a)+4.\end{equation} It
follows that $P_4(x_m)>0$ at $r=r_0$. Consequently, $P_4(x)>0$ on
$[x_m, 1)$ and the proof is now complete. \rule{2mm}{2mm}

\vspace{0.1in}

Recall that steady state $0$ is always locally asymptotically
stable and steady state $1$ is locally asymptotically stable if
$r<r_0$.  Theorems 2.2 and 2.3 provide basins of attractions for
these two attractors when $r<r_0$. Specifically, ${\cal B}(0)=[0,
a)\cup(x_a, \infty)$ and ${\cal B}(1)=(a, x_a)$ whenever $r<r_0$.
Moreover, (2.3) has a period two solution ${\cal S}=\{\bar x_1,
\bar x_2\}$ when $r>r_0$, where $a<\bar x_1<1<\bar x_2$.
 The local stability of ${\cal S}$ can be determined by
the derivative of $f^2(x)$ evaluated at either $\bar x_1$ or $\bar
x_2$. Since $f(\bar x_1)=\bar x_2$, a simple calculation shows
\begin{equation}\begin{array}{l}(f^2(x))^\prime|_{x=\bar x_1}=f^\prime(\bar x_1)f^\prime(\bar
x_2)\\ \\=e^{r((1-\bar x_1)(\bar x_1-a)+(1-\bar x_2)(\bar
x_2-a))}[1+r\bar x_1(1-2\bar x_1+a)][1+r\bar x_2(1-2\bar
x_2+a)].\end{array}\end{equation} \noindent Using (2.11) and
$a<\bar x_1<1$, we have $\bar x_1e^{r(1-\bar x_1)(\bar x_1-a)}<1$
and hence $\bar x_1<x_m$. It follows that  $$1+r\bar x_1(1-2\bar
x_1+a)>0$$ and $$1+r\bar x_2(1-2\bar x_2+a)<1+r\bar
x_2(a-1)<1-r(1-a)<1-r_0(1-a)=-1<0.$$ As a result,
$$(f^2(x))^\prime|_{x=\bar x_1}<0.$$ Therefore, ${\cal S}$ is
locally asymptotically stable if $(f^2(x))^\prime|_{x=\bar
x_1}>-1$ and unstable if $(f^2(x))^\prime|_{x=\bar x_1}<-1$.
Numerically, we can find the  $2$-cycle solution using (2.11) and
then determine its local stability using (2.27) when $r>r_0$. It
is expected that (2.3) undergoes another period-doubling
bifurcation when $(f^2(x))^\prime|_{x=\bar x_1}=-1$. Figure 1(c)
provides a bifurcation diagram for (2.3) when $a=0.5$ and using
$r$ as the bifurcation parameter. The figure conforms with our
analytical findings.

\section{The host-parasitoid model}
\setcounter{equation}{0} In this section, we study the full system
(2.2). Preliminary results are provided in Section 3.1. Stability
and existence of interior steady states when $r<r_0$ and $\beta>1$
are given in Section 3.2. Section 3.3 treats the case for $r>r_0$
and $\beta>1$.
\subsection{Preliminary properties of the host-parasitoid system}
Let $F(x,y)$ denote the map induced by (2.2), i.e.,
\begin{equation} F(x,y)=(f_1(x,y), f_2(x,y)),\end{equation} where
$$f_1(x,y)=xe^{r(1-x)(x-a)-y} \mbox{ and } f_2(x,y)=\beta x(1-e^{-y}).$$
Notice $f_1(x,y)=f(x)e^{-y}$, $x_a>\mbox{max}\{x_m, 1\}$ is the
unique $x$ value for which $f(x_a)=a$, and $x_m$ is the point for
which $f$ attains its maximum. Define the following two regions:
\begin{equation} \begin{array}{l}\Gamma=\{(x,y)\in \mathbb{R}_+^2: x\leq
a \mbox{ and } (x,y)\neq (a, 0)\}\\[1ex]   \Delta=\{(x,y)\in
\mathbb{R}_+^2: x\geq x_a \mbox{ and } (x,y)\neq (x_a,
0)\}.\end{array}\end{equation} We prove that both populations go
extinct if the initial host population is either too small or too
large as given in Proposition 3.1.

\vspace{0.1in}

 \noindent{\bf Proposition 3.1} {\em Solutions of
(2.2) remain nonnegative and are bounded for $t>0$. Moreover,
solutions $(x(t),y(t))$ of (2.2) with $(x(0),
y(0))\in\Gamma\cup\Delta$ converge to $E_0=(0,0)$.}

\vspace{0.1in}

\noindent{\bf Proof.} Let $(x(t), y(t))$ be an arbitrary solution
of (2.2). It is clear that the solution exists and remains
nonnegative for $t>0$. Moreover,  $x(t+1)=f_1(x(t),y(t))\leq
f(x(t))$ for $t\geq 0$. Since $f$ attains its maximum at $x=x_m$,
$x(t+1)\leq f(x_m)$ and thus $y(t+1)\leq\beta x(t)\leq \beta
f(x_m)$ for $t> 0$. Therefore,  solutions are bounded for $t>0$.
If $x(0)\geq x_a$ and $y(0)>0$, then $x(1)=f_1(x(0),y(0))<
f(x(0))\leq f(x_a)=a$ and hence such a solution lies in $\Gamma$
for $t\geq 1$.  Suppose now $x(0)\leq a$ and $y(0)>0$. Then
$x(1)=f_1(x(0),y(0))< f(x(0))\leq x(0)\leq a$ and hence
$\displaystyle\lim_{t\rightarrow\infty}x(t)=0$.  As a result,
$\displaystyle\lim_{t\rightarrow\infty}y(t)=0$ and the solution
converges to $E_0$. \rule{2mm}{2mm}

\vspace{0.1in}

 Let $a<x(0)<1$. Setting
$x(1)=x(0)e^{r(1-x(0))(x(0)-a)}e^{-y(0)}\leq a$ and solving for
$y(0)$, yield
\begin{equation} y(0)\geq
\ln\displaystyle(\frac{x(0)}{a})+(1-x(0))(x(0)-a).
\end{equation}
That is, if $a<x(0)<1$ and $y(0)$ satisfies (3.3), then $x(1)\leq
a$ and the solution converges to $E_0$ by Proposition 3.1. We
summarize below.

\vspace{0.1in}

\noindent{\bf Proposition 3.2} {\em If  $a<x(0)<1$ and $y(0)$
satisfies (3.3), then  the solution converges to $E_0=(0,0)$.}

\vspace{0.1in}

System (2.2) has two more boundary steady states:  $E_1=(a, 0)$
and $E_2=(1,0)$. The Jacobian matrix of (2.2)
evaluated at the boundary steady states are given respectively by \begin{equation}\begin{array}{c}J(E_0)=\left(\begin{array}{cc} e^{-ar} & 0\\
0 & 0
\end{array}
\right), \ J(E_1)=\left(\begin{array}{cc} 1+ar(1-a) & -a\\ 0 &
\beta
a\end{array}\right),\\[3ex] J(E_2)=\left(\begin{array}{cc} 1-r(1-a) & -1\\
0 & \beta
\end{array}\right).\end{array}\end{equation}
Since each of the above Jacobian matrices is triangular, local
stability of these steady states can be easily determined and are
summarized below.

\vspace{0.1in}

\noindent{\bf Proposition 3.3} {\em System (2.2) has three
boundary steady states $E_0=(0,0)$, $E_1=(a, 0)$ and $E_2=(1,0)$,
where $E_0$ is locally asymptotically stable, $E_1$ is a saddle
point if $\beta a<1$ and a repeller if $\beta a>1$. Steady state
$E_2$ is locally asymptotically stable if $r<r_0$ and $\beta<1$
and unstable if either $r>r_0$ or $\beta>1$.}

\vspace{0.1in}

An interior steady state $(x,y)$ of (2.2) satisfies
\begin{equation} y=g(x):=r(1-x)(x-a)
\end{equation}
and \begin{equation} x=h(y):=\displaystyle\frac{y}{\beta
(1-e^{-y})}.\end{equation} It is clear that $g(x)$ is a concave
down parabola with vertex at \begin{equation}\label{xhat}\hat
x=\displaystyle\frac{1+a}{2}\end{equation} and goes through the
points $(a, 0)$ and $(1,0)$. On the other hand,
$$h^\prime(y)=\displaystyle\frac{e^{-y}(e^y-1-y)}{\beta(1-e^{-y})^2}>0
\mbox{ for } y>0$$ and
$$h^{\prime\prime}(y)=\displaystyle\frac{e^{-2y}((y-2)e^{y}+y+2)}{\beta(1-e^{-y})^3}>0
\mbox{ for } y>0.$$ Hence $h(y)$ is increasing and concave up on
$(0, \infty)$ with $\displaystyle\lim_{y\rightarrow
0^+}h(y)=1/\beta$ and $\displaystyle\lim_{y\rightarrow
\infty}h(y)=\infty$. Consequently, if $\beta\leq 1$, then the two
curves $y=g(x)$ and $x=h(y)$ have no positive intersections and
(2.2) has no interior steady state. In such a case,  $\beta
a<\beta\leq 1$ and $E_1=(a,0)$ is a saddle point with a
one-dimensional local stable manifold $\gamma$ by Proposition 3.3.
In the following, we estimate $\gamma$.

Notice that an eigenvector of $J(E_1)$ belonging to the eigenvalue
$\beta a$ can be chosen to be $V=(\displaystyle\frac{a}{1-\beta
a+ar(1-r)}, 1)^T$. Letting $u=x-a$ and $v=y$, system (2.2) is
transformed to
\begin{equation}\label{convsys}
\left \{
\begin{array}{l}
u(t+1)=\left(u(t)+a\right)e^{ru(t)(1-u(t)-a)}e^{-v(t)}-a\\[1ex]
v(t+1)=\beta\left(u(t)+a\right)(1-e^{-v(t)}).
 \end{array}\right.
\end{equation} Steady state $E_1=(a,0)$ of (2.2) becomes steady state
$(0,0)$ of (3.8). Let $v=\gamma(u)=c_1u+c_2u^2+O(u^3)$ for $u$
near $0$ be the local stable manifold of $(0,0)$. Since
eigenvector $V$ is tangent to $\gamma$ at $(0,0)$,  $c_1$ is the
slope of $V$, i.e., $c_1=\displaystyle\frac{1-\beta
a+ar(1-r)}{a}$. Using the positive invariance of $\gamma$ with
respect to system (3.8), we have
\begin{equation*}\begin{array}{l}
v(t+1)=\gamma\left((u(t)+a)e^{ru(t)(1-u(t)-a)-v(t)}-a\right)\\[1ex]
=\beta(u(t)+a)[1-e^{-c_1u(t)-c_2u^2(t)+O(u^3(t))}].
\end{array}\end{equation*}
Equating the expressions with respect to $u$, we  obtain
$c_2=\displaystyle\frac{1+ar(1-a)}{a^2}$, i.e.,
$$\gamma(u)=\displaystyle\frac{1-\beta a+ar(1-a)}{a} u+
\displaystyle\frac{1+ar(1-a)}{a^2}u^2+O(u^3) \mbox{ for } u \mbox{
near } 0.$$ Converting back to the original state variables $x$
and $y$, yields
\begin{equation} \gamma(x)=\beta a-\displaystyle\frac{\beta
a+1+ar(1-a)}{a}x+\displaystyle\frac{1+ar(1-a)}{a^2}x^2+O(x^3)
\end{equation}
for $x$ near $a$.

\vspace{0.1in}

\noindent{\bf Proposition 3.4} {\em Let $\beta\leq 1$. Then system
(2.2) has no interior steady state and  $E_1$ is a saddle point
with the local stable manifold given by (3.9). If $\beta<1$ and
$r<r_0$, then
 $\displaystyle\lim_{t\rightarrow\infty}y(t)=0$ and $\displaystyle\lim_{t\rightarrow\infty}(x(t),y(t))=E_0$, $E_1$, or
$E_2$ for all solutions $(x(t), y(t))$ of (2.2). }

\vspace{0.1in}

\noindent{\bf Proof.} We only need to prove the last statement for
$y(0)>0$. Since $r<r_0$,  solutions of (2.3) satisfy
$\displaystyle\limsup_{t\rightarrow\infty}x(t)\leq 1$ by Theorems
2.2 and 2.3. Hence for any $\epsilon>0$ there exists $t_0>0$ such
that $x(t)<1+\epsilon$ for all $t\geq t_0$ and for all solutions
of (2.3). We can choose $\epsilon>0$ so that $\beta(1+\epsilon)<1$
by the assumption $\beta<1$. Then
$y(t+1)<\beta(1+\epsilon)(1-e^{-y(t)})$ for all $t\geq t_0$. Hence
 $\displaystyle\lim_{t\rightarrow\infty}y(t)=0$ and system (2.2) is asymptotic to the scalar
 equation (2.3). Therefore, solutions of (2.2) either converge to $E_0$, $E_1$,
 or $E_2$ by Theorems 2.2 and 2.3. \rule{2mm}{2mm}

\subsection{Existence and local stability of the interior steady states  when $r<r_0$ and $\beta>1$}

If $r<r_0$, then the asymptotic dynamics of the host population is
completely determined. Moreover, (2.2) has no interior steady
state if $\beta\leq 1$ by Proposition 3.4. Let $\beta>1$. We show
in this subsection that the host-parasitoid system (2.2) has at
most one interior steady state and local stability of the unique
interior steady state is presented. We first prove that (2.2) has
no interior steady state if $0<\displaystyle\frac{1}{\beta}\leq a$
and (2.2) has a unique interior steady state if
$a<\displaystyle\frac{1}{\beta}<1$.

\medskip

\noindent{\bf Theorem 3.5} {\em Let $r<r_0$ and $\beta>1$. System
(2.2) has no interior steady state if
$0<\displaystyle\frac{1}{\beta}\leq a$ and (2.2) has a unique
interior steady state if $a<\displaystyle\frac{1}{\beta}<1$.}

\medskip

\noindent{\bf Proof.} From previous analysis on the isoclines,
(2.2) has at least one interior steady state when
$a<\displaystyle\frac{1}{\beta}<1$. If there exists $\beta>1$ such
that that the two isoclines have two intersections $(x_i, y_i)$,
$i=1, 2$, such that $a< x_1<x_2< 1$, then
$\displaystyle\frac{1}{\beta}$ as a function of $x$
\begin{equation}\label{betaequ}
\displaystyle\frac{1}{\beta}(x)=\displaystyle\frac{x(1-e^{-y})}{y},\end{equation}
where  $y=r(1-x)(x-a)$, has a critical point  in $(x_1, x_2)$.
However, a straightforward calculation shows that
$(\displaystyle\frac{1}{\beta})^\prime(x)>0$ for all $x\in (a, 1)$
if $r\leq r_0$. Indeed,
$$(\displaystyle\frac{1}{\beta})^\prime=\displaystyle\frac{(1-e^{-y})}{y}+x\displaystyle\frac{(1+y)e^{-y}-1}{y^2}\cdot
r(1+a-2x),$$ where $(1+y)e^{-y}-1<0$ for $y>0$ and $rx(1+a-2x)\leq
r(1-a)\leq 2$ for $x\in(a, 1)$ since $r\leq r_0$. Hence
$$(\displaystyle\frac{1}{\beta})^\prime\geq\displaystyle\frac{y-2+(y+2)e^{-y}}{y^2},$$
where $$\left(y-2+(y+2)e^{-y}\right)^\prime=1-(1+y)e^{-y}>0 \mbox{
for } y>0$$ and $$\left(y-2+(y+2)e^{-y}\right)|_{y=0}=0.$$ Thus
$(\displaystyle\frac{1}{\beta})^\prime(x)>0$ for all $x\in (a, 1)$
 and (2.2) has a unique interior steady state if
$a<\displaystyle\frac{1}{\beta}<1$. When $\beta a=1$, then since
$(a,0)$ is itself a steady state, (2.2) has no interior steady
state. Moreover, since $\displaystyle\frac{y}{1-e^{-y}}=\beta x$
is an increasing function of $\beta$ if $x$ is kept fixed, (2.2)
has no interior steady state if $\beta a>1$. \rule{2mm}{2mm}

 \vspace{0.1in}

 The Jacobian matrix of (2.2) evaluated at
any interior steady state $(x,y)$ can be easily shown to be
\begin{equation}J=\left(\begin{array}{cc} 1+xr(1-2x+a) & -x\\
\displaystyle\frac{y}{x} & \beta xe^{-y}
\end{array}
\right),\end{equation} with  \begin{equation} tr
J=1+xr(1-2x+a)+\beta x e^{-y}\end{equation} and
\begin{equation}det J=\beta x
e^{-y}\left(1+xr(1-2x+a)\right)+y,\end{equation} where $y=g(x)$.
Using the second equilibrium equation, we have $
\beta=\displaystyle\frac{y}{x(1-e^{-y})}=\displaystyle\frac{ye^{y}}{x(e^{y}-1)}<\displaystyle\frac{e^{y}}{x}$.
Thus for any interior steady state $(x, y)$ there holds
\begin{equation}\beta x e^{-y}<1.\end{equation} Moreover,
$\mbox{max}\{\displaystyle\frac{1}{\beta}, a\}<x<1$ is always
satisfied.

It is well known  that $(x,y)$ is locally asymptotically stable if
$|tr J|<1+det J$ and $det J<1$ \cite{allen}. Notice $tr J>-1-det
J$ is equivalent to
$$2+xr(1-2x+a)+\beta xe^{-y}(2+xr(1-2x+a))+y>0.$$  Since
$r<r_0$, we have $f^\prime(x)>f^\prime(1)>-1$ if $x\geq x_m$ by
(2.22)  and hence
$$1+e^{r(1-x)(x-a)}(1+xr(1-2x+a))>0.$$ If $x<x_m$, then $1+xr(1-2x+a)>0$ and the
above inequality holds trivially. Therefore,
$2+xr(1-2x+a)>1-e^{-r(1-x)(x-a)}>0$ and $tr J>-1-det J$ is valid
for any interior steady state of (2.2) whenever $r<r_0$. On the
other hand, $tr J<1+det J$ is equivalent to
$$xr(1-2x+a)(1-\beta xe^{-y})<y.$$ This last inequality is
satisfied by (3.14) if $x\geq\hat x$.

The determinant of $J$ given in (3.13) can be rewritten as a
function of $x$:
\begin{equation}
D(x)=\left(1+rx(1+a-2x)\right)\cdot\displaystyle\frac{ye^{-y}}{1-e^{-y}}+y,\end{equation}
where $y=r(1-x)(x-a)$.  Since
$$D(\hat
x)=1\cdot\displaystyle\frac{\hat ye^{-\hat y}}{1-e^{-\hat y}}+\hat
y=\displaystyle\frac{\hat y}{1-e^{-\hat y}}>1, $$
$$D(x)|_{x\geq x_m}\leq 0+y=r(1-x)(x-a)\leq
r(\displaystyle\frac{1-a}{2})^2<0.5 \mbox{ as } r<r_0,$$ $D(x)=1$
has at least one solution on $[\hat x, x_m)$. We prove that $D(x)$
is strictly decreasing on $[\hat x, 1)$. By a direct
differentiation, we have
$$D^\prime(x)=r(1+a-4x)\cdot
\displaystyle\frac{ye^{-y}}{1-e^{-y}}+rx\cdot\displaystyle\frac{d}{dy}(\displaystyle\frac{ye^{-y}}{1-e^{-y}})\cdot
y^\prime+\displaystyle\frac{1-e^{-y}-ye^{-y}}{(1-e^{-y})^2}\cdot
y^\prime$$ $$=r(1+a-4x)\cdot
\displaystyle\frac{ye^{-y}}{1-e^{-y}}+r^2x(1+a-2x)^2\cdot\displaystyle\frac{(1-y)e^{y}-1}{(1-e^{-y})^2}+\displaystyle\frac{1-(1+y)e^{-y}}{(1-e^{-y})^2}\cdot
r(1+a-2x).$$ Since $$1-y<e^{-y} \mbox{ and } e^y>1+y \mbox{ for
all } y>0,$$ $D^\prime(x)<0$ on $[\hat x, 1)$.

If $r\in[r_0/2, r_0)$, then $x_m\leq 1$ and $D(x)=1$ has a unique
positive solution $\underline{x}$ on $(\hat x, x_m)$ with
\begin{equation}\label{D}D(x)\begin{array}{c}>1 \mbox{ on }[\hat x,
\underline{x}) \mbox{ and }
 D(x)<1 \mbox{
on} (\underline{x}, 1).\end{array}\end{equation} If $r\in(0,
r_0/2)$, then $1<x_m$ and we have either $D(x)>1$ on $[\hat x, 1)$
when $\underline x\in(1, x_m)$ or there exists a unique
$\underline{x}\in (\hat x, 1)\subset(\hat x, x_m)$ such that
(\ref{D}) holds.

\noindent{\bf Theorem 3.6} {\em Let $r<r_0$ and
$a<\displaystyle\frac{1}{\beta}<1$ and let $E=(x,y)$ be the unique
interior steady state of (2.2). Then $E$ is locally asymptotically
stable if $x>\underline{x}$ and $E$ is a repeller if
$x<\underline{x}$.}

\medskip

\noindent{\bf Proof.} Assume first $x>\underline{x}$. We can
exclude the discussion for $\underline x\in(1, x_m)$ since $(x,y)$
is not an interior steady state for $x>\underline x>1$. Then $det
J<1$ by (3.16) and $tr J>-1-det J$ is satisfied since $r<r_0$ by
the previous discussion. Furthermore, $tr J<1+det J$ is also valid
as $x>\underline x>\hat x$. Therefore $|tr J|<1+det J$ is proven
and $E$ is locally asymptotically stable if $x>\underline{x}$.

Suppose now $x<\underline{x}$. We separate the discussion into
$x\leq\hat x$ and $x>\hat x$. If $x\leq \hat x$, then since
$x>\displaystyle\frac{1}{\beta}$, (3.13) implies
$$det J=\beta x+r\beta x^2e^{-y}(1-2x+a)\geq\beta x>1.$$   To prove $E$ is a repeller, it is
equivalent to verify that both of the eigenvalues $\lambda_\pm$ of
$J$ have modulus greater than $1$.  If $(tr J)^2-4 det J<0$, then
the eigenvalues  are complex with  modulus greater than $1$ by the
fact that $|\lambda_{\pm}|^2=|\lambda_+\lambda_-|=detJ>1$. If $(tr
J)^2-4 det J\geq 0$, then both  $\lambda_{\pm}$ are real and $(tr
J)^2\geq 4det J>4$. In this case, $\mbox{min}|\lambda_\pm|>1$ if
and only if $|tr J|<1+ det J$. Since $x\leq \hat x<x_m$, we have
$tr J>0$. Therefore,  $|tr J|=tr J<1+ det J$ is equivalent to
\begin{equation} r\left(x(1+a-2x)-(1-x)(x-a)-x(1+a-2x)\cdot\displaystyle\frac{y}{e^y-1}\right)<0,\end{equation}
where $y=r(1-x)(x-a)$.

Let the left hand side of the above inequality be denoted by
$V(x)$, $a<x\leq\hat x$. Then $\displaystyle\lim_{x\rightarrow
a^+}V(x)=0$ since $\displaystyle\lim_{y\rightarrow
0^+}\displaystyle\frac{y}{e^y-1}=1$. If we can show that
$V^\prime(x)<0$ on $(a, \hat x)$, then $V(x)<0$ on $(a, \hat x)$
is proved and $E$ is a repeller. A straightforward differentiation
yields
$$V^\prime(x)=-r\left((1+a-2x)\displaystyle\frac{y}{e^y-1}\left(1-\displaystyle\frac{xr(1-2x+a)}{2}\right)+B\right),$$
where
$$B=x\left(2(1-\displaystyle\frac{y}{e^y-1})+r(1-2x+a)^2\cdot\left(\displaystyle\frac{(2-y)e^y-(2+y)}{2(e^y-1)^2}\right)\right).$$
Since
$$\left((2-y)e^y-(2+y)\right)^\prime=e^y\left((1-y)-e^{-y}\right)<0
\mbox{ for } y>0$$ and $$\left((2-y)e^y-(2+y)\right)|_{y=0}=0,$$
we have $(2-y)e^y-(2+y)<0$ for $y>0$. Also,
$r(1-2x+a)^2<r_0(1-a)^2<2$ on $(a, \hat x]$. Hence for $a<x\leq
\hat x$,
$$B\geq
2x\left((1-\displaystyle\frac{y}{e^y-1})+\displaystyle\frac{(2-y)e^y-(2+y)}{2(e^y-1)^2}\right)
=2x\left(\displaystyle\frac{2e^{2y}-(2+3y)e^y+y}{2(e^y-1)^2}\right).$$
Since
$$(2e^{2y}-(2+3y)e^y+y)^\prime=4e^{2y}-(5+3y)e^y+1$$
$$\geq\left(4(1+y)-(5+3y)\right)e^y+1>e^y\left(e^{-y}-(1-y)\right)>0$$
for $y>0$ and $(2e^{2y}-(2+3y)e^y+y)|_{y=0}=0$, we have
 $2e^{2y}-(2+3y)e^y+y>0$ for $y>0$, i.e., for $x\in(a, \hat x]$.
 This shows that $B>0$ on $(a, \hat x]$. Moreover, $1+a-2x<1-a$ and hence
$\displaystyle\frac{xr(1-2x+a)}{2}<\displaystyle\frac{r(1-a)}{2}<1$
for all $x\in (a, \hat x]$.
  Therefore, $V^\prime(x)<0$ on $(a, \hat x]$ and $V(x)<0$ on $(a, \hat x]$ is verified. It follows that the interior steady
state $E$ is a repeller if $x\leq \hat x$. The proof of $x\in(\hat
x, \underline{x})$ is trivial since $det J>1$ by (3.16) and   $tr
J>-1-det J$ always holds. Moreover, $tr J<1+det J$ from an earlier
observation since  $x>\hat x$. Therefore, $|tr J|<1+det J$ and $E$
is a repeller if $x\in(\hat x, \underline{x})$.  \rule{2mm}{2mm}

\medskip

 In particular, the interior steady state $E=(x,y)$ is locally asymptotically stable when $x$ is near $1$.
 Finally, we prove that steady state $E_0=(0,0)$ is
globally asymptotically stable in $\{(x,y)\in\mathbb{R}_+^2:
y>0\}$ if $\beta a>e^{\frac{1-a}{2}}$ and $r<r_0$.

\medskip

\noindent{\bf Theorem 3.7} {\em If $\beta a>e^{\frac{1-a}{2}}$ and
$r<r_0$, then $E_0=(0,0)$ is globally asymptotically stable in
$\{(x,y)\in\mathbb{R}_+^2: y>0\}$ for system (2.2).}

\medskip

\noindent{\bf Proof.} Let $(x(t), y(t))$ be an arbitrary solution
of (2.2) with $x(0), y(0)>0$. Then $x(t), y(t)>0$ for $t\geq 0$.
If there exists $t_0>0$ such that $x(t_0)\leq a$, then
$\displaystyle\lim_{t\rightarrow\infty}(x(t),y(t))=E_0$ by
Proposition 3.1. Assume now $x(t)>a$ for all $t\geq 0$. Since
$r(1-x)(x-a)<r_0(1-\hat x)(\hat x-a)<\displaystyle\frac{1-a}{2}$
for all $ x\geq a$, we have $x(t+1)\leq
x(t)e^{\frac{1-a}{2}}e^{-y(t)}$ and $y(t+1)\geq \beta
a(1-e^{-y(t)})$ for all $t>0$. Consider
\begin{equation}\label{scalar} z(t+1)= \beta a(1-e^{-z(t)}).\end{equation}
The scalar equation (\ref{scalar}) has a unique interior steady
state $\bar z$ and $z(t)\rightarrow \bar z$ as $t\rightarrow
\infty$ if $z(0)>0$. It follows that
$\displaystyle\liminf_{t\rightarrow\infty}y(t)\geq \bar z$.

 We verify that $\bar z>\displaystyle\frac{1-a}{2}$. Let $S(z)=\beta
a(1-e^{-z})-z, \ z\geq 0$. Then $S(0)=0$ and $S^\prime(z)=\beta
ae^{-z}-1>0$ for $0\leq z\leq \displaystyle\frac{1-a}{2}$ by the
assumption, i.e., $S(z)>0$ for
$0<z\leq\displaystyle\frac{1-a}{2}$. Hence $\beta a (1-e^{-z})>z$
for $0<z\leq\displaystyle\frac{1-a}{2}$ and $\bar
z>\displaystyle\frac{1-a}{2}$ is shown. For any $\epsilon>0$ there
exists $t_1>0$ such that $y(t)>\bar z-\epsilon$ for all $t\geq
t_1$. We can choose $\epsilon>0$ such that
$\displaystyle\frac{1-a}{2}+\epsilon-\bar z<0$. Then $$x(t+1)\leq
x(t)e^{\frac{1-a}{2}+\epsilon-\bar z} \mbox{ for } t\geq t_1$$
implies
 $\displaystyle\lim_{t\rightarrow\infty}x(t)=0$. Therefore, the
solution converges to $E_0$ and $E_0$ is globally asymptotically
stable in $\{(x,y)\in\mathbb{R}_+^2: y>0\}$ since the solution was
arbitrary and $E_0$ is locally asymptotically stable.
\rule{2mm}{2mm}

\medskip

Dynamical behavior of system (2.2) is now summarized for $r<r_0$
and is illustrated in Figure 2. We do not  analytically prove
local bifurcations of (2.2).
\begin{itemize}
\item $\beta<1$: (2.2) has no interior steady state and solutions
converge to one of the boundary steady states $E_0$, $E_1$ or
$E_2$ depending on initial conditions.  It is expected that the
system undergoes a transcritical bifurcation at $\beta=1$ and
(2.2) has an asymptotically stable interior steady state when
$\beta>1$ and is near $1$.

\item $1<\beta<\displaystyle\frac{1}{a}$: (2.2) has a unique
interior steady state $E=(x,y)$, where $E$ is locally
asymptotically stable if $x>\underline x$ and is a repeller if
$x<\underline x$. As $|tr J|<1+det J$ holds for any interior
steady state $E=(x,y)$ while $D(\underline x)=1$ and
$D^\prime(x)<0$ on $[\hat x, 1)$, one expects that a
Neimark-Sacker bifurcation  occurs when $x=\underline x$. Recall
that $\displaystyle\frac{1}{\beta}$ is a strictly increasing
function of $x$. Therefore, there exists a unique $\beta_c$,
$1<\beta_c<\displaystyle\frac{1}{a}$, that corresponds to the
unique $\underline x$. Moreover, $E$ is locally asymptotically
stable if $1<\beta<\beta_c$ and $E$ is a repeller if
$\beta_c<\beta<\displaystyle\frac{1}{a}$. A Neimark-Sacker
bifurcation occurs at $\beta=\beta_c$ and the system has
quasi-periodic solutions when $\beta>\beta_c$ and is close to
$\beta_c$. In the numerical example given in Figure 2 for $a=0.5$
and $r=2.0$, we have $\beta_c\approx 1.3777$ and (2.2) has locally
stable quasi-period solutions when $\beta\in(\beta_c, 1.42812)$.
We expect that (2.2) undergoes a saddle node bifurcation when
$\beta a=1$ so that the interior steady state $E$ disappears as
$\beta$ is increased beyond $1/a$.

\item $\beta>1/a$: (2.2) has no interior steady state and  $E_0$
is globally asymptotically stable if $\beta a>e^{\frac{1-a}{2}}$.
Numerical simulations suggest that this global behavior remains
true  when $\beta a>1$.
\end{itemize}

{\bf Remark.} From the above summary it is known  that for a fixed
$r$, $0<r<r_0$, there exists a unique $\beta_c\in (1, 1/a)$ such
that $E$ is locally asymptotically stable if $\beta\in(1,
\beta_c)$ and $E$ is a repeller when $\beta\in(\beta_c, 1/a)$. A
Neimark-Sacker bifurcation occurs at $\beta=\beta_c$ and the
system has quasi-periodic solutions when $\beta>\beta_c$ and is
close to $\beta_c$. As $r$ increases, the critical value $\beta_c$
also increases according to the numerical investigations and
therefore the $\beta$ parameter region for which $E$ is locally
asymptotically stable becomes larger. For instance, when $a=0.5$,
we have $\beta_c\approx 1.318$, $1.3378$, $1.35765$,  $1.3777$,
$1.41819$,  $1.45519$, $1.4589$ as $r$ increases from $0.5$,
$1.0$, $1.5$, $2.0$, $3.0$, $3.9$ to $3.99$ respectively. So if we
fix $\beta=1.5$ and $r=3.99$, then since $\beta$ is close to
$1.4589$, the system has a  stable quasi-periodic solution as
shown in Figure 3(a) with initial condition $(0.7512, 0.2437)$
chosen to be close to the interior steady state. On the other
hand, if $r$ is decreased to $r=0.5$, then since $\beta=1.5$ is
somewhat much larger than the $\beta_c$ value of $1.318$, dynamics
of (2.2) behave much like Figure 2(d). Solutions spiral away from
the interior steady state and move toward $E_0$. We also plot the
stable invariant loop solution in Figure 3(b) for $r=0.5$,
$\beta=1.325$ and initial condition $(0.7811, 0.0308)$.

\subsection{Existence and local stability of the interior steady states  when $r>r_0$ and $\beta>1$}
In this subsection we study stability of interior steady states
using  Jury conditions \cite{allen}. Let $E=(x,y)$ be an interior
steady state of (2.2). Recall that such an interior steady state
exists if $a<1/\beta<1$. The Jacobian matrix $J(E)$ evaluated at
$E$ is given by (3.11) with trace and determinant given by (3.12)
and (3.13) respectively. We separate the discussion into two
cases: $a<x<\hat x$ and $\hat x\leq x<1$.

Assume first $\hat x\leq x<1$ and let $D(x)$ be the determinant of
$J(E)$ given by (3.15). Then $D(\hat x)=\displaystyle\frac{\hat
y}{1-e^{-\hat y}}>1$  and $\displaystyle\lim_{x\rightarrow
1^-}D(x)=1-r(1-a)<1-r_0(1-a)=-1.$ By the same argument as in the
case of $r<r_0$ and $\beta>1$, it can be verified that
$D^\prime(x)<0 \mbox{ for } x\in [\hat x, 1).$ Hence there exists
a unique $\check{x}$ in $(\hat x, 1)$ such that $D(\check x)=1$
and
\begin{equation}D(x)>1 \mbox{ on }[\hat x,
\check{x})\mbox{ and } D(x) <1 \mbox{ on } (\check{x},
1].\end{equation} Moreover, since $x\geq \hat x$, $tr J<1+det J$
holds trivially.

Let \begin{equation} T(x)=tr J+1+det
J=\left(2+rx(1+a-2x)\right)(1+\displaystyle\frac{y}{e^y-1})+y.\end{equation}
Then $-1-det J<tr J$ is equivalent to $T(x)>0$. A direct
differentiation yields
$$T^\prime(x)=r(1+a-2x)\cdot
A-2x(1+\displaystyle\frac{y}{e^y-1})+r^2x(1+a-2x)^2\cdot\displaystyle\frac{(1-y)e^y-1}{(e^y-1)^2},$$
where
$$A=2+\displaystyle\frac{y}{e^y-1}+2\cdot\displaystyle\frac{(1-y)e^y-1}{(e^y-1)^2}.$$
Observe that $(1-y)e^y-1<0 \mbox{ for } y>0$, and
$$ A>\displaystyle\frac{1}{(e^y-1)^2}\left[(2(1+y)-(y+2))e^y-y\right]
=\displaystyle\frac{y(e^y-1)} {(e^y-1)^2}>  \mbox{ for }  y>0.$$
Therefore $T^\prime(x)<0 \mbox{ for } x\in[\hat x, 1).$ Moreover,
$T(\hat x)=2(1+\displaystyle\frac{\hat y}{e^{\hat y}-1})+\hat y>0$
and $\displaystyle\lim_{x\rightarrow 1^-}T(x)<0\mbox{ if } r>r_0.$
Hence there exists a unique $\tilde x\in(\hat x, 1]$ such that
\begin{equation}T(x)>0\mbox{ on } [\hat x, \tilde x)\mbox{ and }
T(x) <0 \mbox{ on } (\tilde x, 1].\end{equation} In general, we
assume $\check x\neq \tilde x$. Moreover, we do not consider the
critical cases when either $x=\check x$ or $x=\tilde x$.

\vspace{0.1in}

\noindent{\bf Theorem 3.8} {\em Let $r>r_0$, $\beta>1$ and let
$E=(x,y)$ be an interior steady state of (2.2) with $\hat x\leq
x<1$. The following statements hold for $E$.
\begin{itemize}\item[(a)] If $\tilde x<\check x$, then (i) $E$ is a repeller if $\hat x\leq x<\tilde x$, (ii) $E$
is a saddle point if $\tilde x<x<\check x$, and (iii)  $E$ is a
saddle point if $x>\check x$. \item[(b)] If $\check x<\tilde x$,
then (i) $E$ is a repeller if $\hat x\leq x<\check x$, (ii) $E$ is
locally asymptotically stable if $\check x<x<\tilde x$, and (iii)
$E$ is a saddle point if $x>\tilde x$.\end{itemize}}

\vspace{0.1in}

\noindent{\bf Proof.}  (a) If $\hat x\leq x<\tilde x$, then
$T(x)>0$ by (3.21) and $D(x)>1$ by (3.19). Therefore, $E$ is a
repeller by a similar proof as that of Theorem 3.6. This proves
(i). If $\tilde x<x<\check x$, then $T(x)<0$ and $D(x)>1$.
Therefore, $tr J<-2$, $(tr J)^2-4det J>(1+det J)^2-4 det J=(1- det
J)^2>0$ and eigenvalues $\lambda_{\pm}$ of $J$ are real numbers.
Since $T(x)<0$ and $det J>1$, a direct computation shows
$\lambda_+<1$ and $\lambda_-<-1$. Moreover,  $\lambda_+>-1$ is
satisfied since $-2-tr J>0$ and $T(x)<0$. Therefore, $E$ is a
saddle point, which proves (ii). Suppose now $x>\check x$. Then
$D(x)<1$ and $T(x)<0$. If $-1-det J\leq 0$, then $(tr J)^2-4 det
J>0$. If $-1-det J>0$, then $det J<0$. Hence  eigenvalues
$\lambda_{\pm}$ of $J$ are real numbers. It is then
straightforward to show that $E$ is a saddle point.

The proof of (b) is similar to (a). If inequalities in (i) are
satisfied, then $D(x)>1$, $T(x)>0$ and thus $E$ is a repeller by
an argument similar to the proof of Theorem 3.6. If $x\in (\check
x, \tilde x)$, then $D(x)<1$ and $T(x)>0$. Therefore $E$ is
locally asymptotically stable by an earlier observation. If $x\in
(\tilde x, 1)$, then $D(x)<1$ and $T(x)<0$. Therefore, $E$ is a
saddle point by the proof of (a)(iii). \rule{2mm}{2mm}

\vspace{0.1in}

If the $x$ component of an interior steady state is such that
$x<\hat x$, then $$det J=\left(1+rx(1+a-2x)\right)\beta
xe^{-y}+y>\beta x e^{-y}+y=\displaystyle\frac{y}{1-e^{-y}}>1$$ and
such an interior steady state is always unstable.

\vspace{0.1in}

\noindent{\bf Proposition 3.9} {\em Let $r>r_0$, $\beta>1$ and let
$E=(x,y)$ be an interior steady state of (2.2) with $x< \hat x$.
Then $E$ is unstable. }

\vspace{0.1in}

However, since $r>r_0$, the function
$\displaystyle\frac{1}{\beta}(x)$ defined in (\ref{betaequ}) may
not be strictly increasing on $(a, 1)$ as in the case of $r<r_0$.
Therefore, (2.2) may have multiple interior steady states.
Furthermore, let \begin{equation}\label{Y}
Y(x)=g(x)-h^{-1}(x),\end{equation} where $g$ and $h$ are defined
in (3.5) and (3.6) respectively. Then a tedious  computation shows
that \begin{equation}\label{3DY} Y^{\prime\prime\prime}(x)>0
\mbox{ on } (0, \infty).\end{equation} Therefore, $Y(x)=0$ has at
most three positive solutions. It follows that (2.2) has at most
three interior steady states. There are plenty of parameter values
for which there are two interior steady states. Although the
analytical result indicates that the system can have at most three
interior steady states, we do not find the existence of three
interior steady states numerically. We conjecture that the system
has at most two interior steady states when $r>r_0$ and
$1<\beta<1/a$.

Dynamics of (2.2) are complicated as illustrated by a numerical
example with $a=0.1$ and $r=3.0>r_0$ given in Figure 3. Notice
that both populations become extinct  when $\beta$ is  large. It
is suspected that the large turnover of parasitoids will drive the
host population to below its Allee threshold so that the host
population goes extinct and then the parasitoid population will
inevitably become extinct.

 The following table provides
notation used in Sections 2 and 3 and Table 2 summarizes existence
and stability of the interior steady states.

\vspace{0.1in}

\begin{center} Table 1: Notation \end{center}

 \begin{tabular}{|c|l|}\hline
 Notation & Definition\\[1ex]
 \hline
 $\hat x$ & $\displaystyle\frac{1+a}{2}$\\[1ex]
 \hline
 $x_m$ & The critical point of $f(x)$\\[1ex]
 \hline
 $x_a$ & The $x$ value for
 which $f(x)=a$, $x_a>\mbox{max}\{x_m ,1\}$\\[1ex]
 \hline
 $\underline x$ & The $x$ value for which $D(x)=1$ when $r<r_0$ and
 $\beta>1$\\[1ex]
 \hline
$\check x$ & The $x$ value for which $D(x)=1$ when $r>r_0$ and
 $\beta>1$\\[1ex]
 \hline
 $\tilde x$ & The $x$ value for which $T(x)=0$ when $r>r_0$ and
 $\beta>1$\\[1ex]
 \hline
 \end{tabular}

\vspace{0.1in}

\begin{center} Table 2: Existence and stability of an interior steady state $E=(x,y)$ \end{center}
 \begin{tabular}{|c|l|l|}\hline
$r$ & Existence & Stability\\[1ex]
\hline $r<r_0$ & $a<\displaystyle\frac{1}{\beta}<1$ &
\begin{tabular}{l} Locally asymptotically stable
if $x>\underline x$\\[1ex]
 Repeller  if $x<\underline x$\end{tabular}\\[1ex] \hline
 $r>r_0$ & $a<\displaystyle\frac{1}{\beta}<1$ & \begin{tabular}{l}Locally asymptotically stable if $\check x<x<\tilde x$ \\[1ex]
 Saddle point if $\tilde x<x<\check x$ or $x>\check
 x>\tilde x$ or $x>\tilde x>\check x$\\[1ex]
Repeller
 if $\hat x\leq x<\tilde x<\check x$ or $\hat x\leq x<\check
 x<\tilde x$\\[1ex]\end{tabular}\\
 \hline
\end{tabular}

\section{Discussion}
\setcounter{equation}{0} Hosts and parasitoids  are frequently
insect populations with distinctive life stages and therefore
discrete-time models are appropriate to describe such populations.
Allee effects are biological phenomena which arise when reduced
fitness or declined population growth occurs at low population
densities or sizes. Mechanisms of Allee effects include failure of
finding mates to reproduce or lack of cooperative individuals to
explore resources collaboratively. When a population is subject to
Allee effects, it is well known that there exists a critical
population level  below which the population will go extinct
\cite{courchamp}. Consequently, studying population interactions
involving
 Allee effects becomes an important subject of contemporary
 research.

 In this investigation, we first study a single host population that
 follows the classical Ricker function and is also subject to Allee effects.
The incorporation of Allee effects into the Ricker equation is
different than those of the Ricker models with Allee effects in
\cite{elaydi1, jang,  kang, li,  liva} but is similar to that
given in \cite{allen1, liebhold}. It is proved that there exists a
threshold $r_0$ in terms of the Allee threshold $a$ such that the
population either goes to extinction or stabilizes at the carrying
capacity if its intrinsic growth rate $r$ is less than $r_0$. The
basins of attraction of the local attractors are determined in
Theorems 2.2 and 2.3. A period-doubling bifurcation occurs at
$r=r_0$ (cf. Proposition 2.1), and, by numerical simulations, the
host population will undergo a cascade of period-doubling
bifurcations and eventually be chaotic if $r$ is increased
further.

For the host-parasotoid system (2.2), it is shown in Proposition
3.1 that both populations will go extinct if the initial host
population is either less than $a$ or is greater than the level of
$x_a$. The only interesting dynamical behavior occurs when the
initial host population is between $a$ and $x_a$. Moreover, the
parasitoids will become extinct if $\beta$ is less than $1$ (cf.
Proposition 3.4). When $\beta$ exceeds $1$ and $r<r_0$, then the
interaction can support at most one interior steady state by
Theorem 3.5. The interior steady state can be either locally
asymptotically stable or is a repeller  as illustrated in Theorem
3.6. If $r>r_0$, then the system may have multiple interior steady
states. We provide linear stability analysis of these interior
steady states in Theorem 3.8 and Proposition 3.9.

To compare the present investigation with other studies mentioned
in the Introduction, we notice that \cite{dennis} is on stochastic
models while \cite{cushing, zhou, zhou1} deal with continuous-time
models. Discrete-time models of Allee effects include
\cite{elaydi1, jang1, jang, kang, li, liva, schreiber}, where
\cite{elaydi1, li, schreiber} are of single population models and
\cite{kang, liva} are of competition models. Moreover, the
predator-prey model in \cite{jang1} uses a Beverton-Holt growth
equation. We therefore compare the present study with the model in
\cite{jang}:
\begin{equation}\label{compare}\left\{\begin{array}{l}
N(t+1)=N(t)e^{r-N(t)}\displaystyle\frac{N(t)}{m+N(t)}e^{-bP(t)}\\[2ex]
P(t+1)=\beta N(t)(1-e^{-bP(t)}),
\end{array}\right.
\end{equation}
where $1/m>0$ is an individual host's searching efficiency. There
exists a threshold $\tilde r_0$ as an increasing function of $m$,
such that $E_0=(0,0)$ is globally asymptotically stable if
$r<\tilde r_0$.  The system  has two more boundary states
$E_{1i}=(\bar N_i, 0)$, $i=1, 2$, when $r>\tilde r_0$, where $\bar
N_1<\bar N_2$ depends on both $r$ and $m$. The host population
goes extinct and so does the predator population if the initial
host population is smaller than $\bar N_1$. System (4.1) may have
multiple interior steady states and the system  can either undergo
a period-doubling bifurcation or a Neimark-Sacker bifurcation when
an interior steady state loses its stability. Observe that global
extinction also occurs in system (2.2) when $r<r_0$ and $\beta
a>1$. The extinction of both populations in (2.2) is likely due to
the predator's overexploitation as $\beta$ has to be large. Global
asymptotic stability of $E_0$ in (4.1), however,  is more likely
to occur if an individual host's searching efficiency is small,
i.e., if $m$ is large, so that $r<\tilde r_0$ is more likely to
hold. Furthermore, (2.2) can only undergo a Neimark-Sacker
bifurcation when $r<r_0$ while (4.1) can exhibit a period-doubling
bifurcation (cf. Fig. 3(c)--(d) of [11]) when an interior steady
state loses its stability.

 \vspace{0.1in}

\noindent{\bf Acknowledgements.} We thank both referees for their
valuable comments and suggestions. Y. Chow acknowledges the
support of the National Science Council of Taiwan. S. Jang thanks
the Institute of Mathematics, Academia Sinica, Taiwan for its
financial and staff support for her winter 2012 and summer 2013
visits.

\end{document}